\documentclass[a4paper,14pt]{article}
\usepackage{amsmath, amsthm, amsfonts, amssymb}

\theoremstyle{plain}
\newtheorem{thm}{Theorem}
\newtheorem*{thm*}{Theorem}
\newtheorem{lem}{Lemma}
\theoremstyle{definition}
\newtheorem{defn}{Definition}

\theoremstyle{remark}
\newtheorem{remk}{Remark}
 
\newcommand{\1}{1\!\!\,{\rm I}}

\newcommand{\wt}{\widetilde}

\newcommand{\ov}{\overline}

\newcommand{\ve}{\varepsilon}

\newcommand{\mbR}{{\mathbb R}}
\newcommand{\mbN}{{\mathbb N}}

\newcommand{\Z}{\mathbb{Z}}

\newcommand{\Ce}{\mathrm{C}([0,T])}

\newcommand{\sign}{\mathop{\mathrm{sign}}}

\newcommand{\cF}{{\cal F}}

\newcommand{\cB}{{\cal B}}

\newcommand{\supp}{\supp{\rm supp}}

\newcommand{\E}{\mathrm{E}}
\newcommand{\Pb}{\mathrm{P}}

\newcommand{\de}{\,\mathrm{d}}
\newcommand{\Sum}{\sum\limits}

\newcommand{\Lim}{\lim\limits}

\renewcommand{\le}{\leqslant}
\renewcommand{\ge}{\geqslant}
\newcommand{\peq}{\stackrel{\mathrm{d}}{=}}

\newcommand{\mbZ}{\mathbb{Z}}

\begin{document}

\large

\title{ On a limit behavior of a sequence of Markov processes perturbed in a neighborhood of a singular point}
\author{Andrey Pilipenko\footnote{Institute of Mathematics, National
Academy of Sciences of Ukraine, Kiev, Ukraine; The National Technical University of Ukraine ``Kyiv Polytechnic Institute'', Kiev, Ukraine
\newline e-mail: pilipenko.ay@yandex.ua}
  \ \ and \ \
  Yuriy Prykhodko  \footnote{The National Technical University of Ukraine ``Kyiv Polytechnic Institute'', Kiev, Ukraine
\newline e-mail: prykhodko@matan.kpi.ua}
 }

\maketitle

We study a limit behavior of a sequence of Markov processes (or Markov chains)  such that their distributions outside of any neighborhood of a ``singular'' point attract to some probability law. In any neighborhood of this point the behavior may be irregular. As an example of the general result we consider a symmetric random walk with the unit jump that is perturbed in a neighborhood of 0. The invariance principle is obtained under standard scaling of time and space. The limit process turns out to be  a skew Brownian motion.

\section{Introduction}
The main goal is to investigate a limit behavior of a sequence of Markov processes   $\{X_n\}$ (or scaled Markov chains).We assume  that outside any neighborhood of a fixed  ``singular'' point the distribution of the sequence converges to a certain law. At a neighborhood of this singular point the  behavior of a sequence may be irregular. For example,
let  $S_n=\xi_1+\dots+\xi_n, n\ge 0$ be a standard random walk with $E\xi_k=0$ and  finite variance of the jump $\sigma^2=Var\xi_k.$ It is well known that the distributions of  $\{\frac{S_{[nt]}}{\sigma\sqrt{n}},\ n\ge 1\}$ converge weakly to a Brownian motion in  $D([0;T])$. Assume now that   $\{\tilde S_n\}$ is a Markov chain such that its transition probabilities outside the set 
$[-m,m] $ coincide with transition probabilities of  $\{S_n\},$  and they are arbitrary inside of $[-m,m] $. Then a sequence $\{\frac{\tilde S_{[nt]}}{\sigma\sqrt{n}},\ n\ge 1\}$ behaves exactly as $\{\frac{S_{[nt]}}{\sigma\sqrt{n}},\ n\ge 1\}$ outside any neighborhood of zero (note that the rescaled interval  $[-m/\sqrt{n}, m/\sqrt{n}]$ is shrinking to zero as   $n\to\infty$). In this case a limit of $\{\frac{S_{[nt]}}{\sigma\sqrt{n}},\ n\ge 1\}$ may be not a Brownian motion. For example, if   $\{\tilde S_n\}$ is a walk on $\mbZ$ with the following transition probabilities  $\tilde p_{i,i\pm 1}=1/2$ for $i\neq0$, and $ p_{0,1}=p, p_{0,-1}=q=1-p$, then the limit is a skew Brownian motion  \cite{HarrShepp81}. That is, the limit is a continuous Markov process with transition probability density function 
\begin{equation*}
p_t(x,y) = \varphi_t (x-y) + \gamma\sign(y)\varphi_t (|x|+|y|),\
x,y \in \mathbb{R},
\end{equation*}
where $\varphi_t(x) = \frac{1}{\sqrt{2\pi t}} e^{-x^2/2t}$ is the density of the normal distribution  $N(0,t)$. In this case parameter of skewness  $\gamma$  equals  $(p-q).$ The case of the arbitrary  $m$ and bounded integer-valued  jumps was considered in \cite{Minlos97, Yarotskij, PP11}. The limit process also was a skew Brownian motion.  If jumps from $[-m,m]$ are not integrable, then the limit could be a discontinuous process, for example a Brownian motion with Wentzell's boundary conditions \cite{PP2014}.

Another examples   are diffusions on graphs, see  \cite{EnriquezKifer, WentzelFreidlin}. They can be obtained, for example, as a limit
of a symmetric random walks on edges of a graph with additional assumptions on jumps from vertices. Diffusions on graphs can also be obtained as a limit of more complex systems, where in a limit few vertices are   ``tied'' together \cite{Kulik}.

Construction of  a reflecting Brownian motion in a cone is also can be considered as a diffusion with irregularity at a vertex, see \cite{VW, Kwon, KwonWilliams} and references therein. Until the hitting of the vertex, reflecting Brownian motion can be defined as a strong solution of a certain SDE. A problem to determine a process after the hitting of the vertex is non-trivial. One way of a construction of the  process is to determine a limit   of a sequence of reflecting Brownian motions that   jumps instantly on a vector with length $\ve$ inside the cone when these processes hit  the vertex. Certainly it is not obvious that the limit exists and is unique.

One more example  is a Peano phenomenon for diffusions with a small noise, see for example  \cite{BaficoBaldi, KrykunMakhno}. There one studied a limit behavior as $\ve\to 0$  of the following   diffusions 
$$
d\xi_\ve(t)=b(\xi_\ve(t))dt+\ve\sigma(\xi_\ve(t))dw(t).
$$
It was assumed  that   $b$ is continuous,   $b(0)=0$, $\xi(0)=0,$ but the limit equation 
\begin{equation}\label{eqKrykun}
d\xi(t)=b(\xi(t))dt
\end{equation}
may have non unique solution. Under some assumptions on  $b$ and $ \sigma$ it was proved that the limit process equals the maximal solution of \eqref{eqKrykun}
with some probability $p$, and equals the minimal solution with probability  $(1-p)$.

We propose a technique to determine  the limit of   processes with irregularities at a neighborhood of a fixed point. Generally, we even do not assume that the initial sequence   $\{X_n\} $ is Markov. A general abstract result is given in   \S 2. We provide sufficient conditions if $\{X_n\} $ satisfies some renewal properties in   \S 3 (Theorem \ref{prop1}). In \S 4 we give an example of our methodology application. There a Markov chain $\{\tilde S_k, k\geq 0\}$ in $\mbZ$ is considered, where transition probabilities are
$$
p_{i,i\pm 1}=1/2,\ \ \mbox{for} \ \ \  |i|>m,
$$
$$
\sum_jp_{ij}|j|<\infty,\ |i|\le m,
$$
i.e., if $|i|\le m$, then the expectation of a jump is finite.

It is proved that a 
sequence $\{\frac{\tilde S_{[nt]}}{\sigma\sqrt{n}},\ n\ge 1\}$ converges to a skew Brownian motion. A parameter of skewness is calculated exactly. In contrast to  \cite{Minlos97, Yarotskij, PP11} we do not assume the boundedness of jumps for $|i|\le m.$
This example may be considered by resolvent or semigroup theory. However, our considerations have more transparent probabilistic interpretation.  On the other hand, our methods can be applied also for continuous semi-Markov processes (see the corresponding definition in  \cite{Harlamov}) that are perturbed at   neighborhoods of a point  $x^*$. We only need to assume that  some exit moments and entrance moments are renewal.

\section{A general limit result
}

In this section we describe a general approach for   studying of a limit behavior of  stochastic processes with irregularities at a neighborhood of a fixed point $x^*$. In order to determine a limit we will  ``delete'' certain parts of trajectories that belong to a neighborhood of $x^*$ and then study a sequence with  ``deleted'' time intervals.

Consider the formal reasoning.

 Let $(E,\rho)$ be  a locally compact metric space. By $C([0,\infty))=C_E([0,\infty))$ and $D([0,\infty))=D_E([0,\infty))$ denote spaces of functions with values in   $E$, that have continuous and c\`{a}dl\`{a}g trajectories respectively.

Let $f\in D([0, \infty))$,
$\sigma=\{\sigma_n\},\ \tau=\{\tau_n\}$ be sequences of real numbers such that
\begin{equation}\label{eq1.1}
0\le\tau_0\le\sigma_0<\tau_1<\sigma_1<\tau_2<\ldots,
\end{equation}
\begin{equation}\label{eq1.2}
\lim_{n\to\infty}\tau_n=+\infty, \ \lambda(\cup_k[\sigma_k, \tau_{k+1}))=+\infty,
\end{equation}
where  $\lambda$ is the Lebesque measure.

Define
\[
L(t)=L_{\tau, \sigma}(t):=\int^t_0\1_{\cup_k[\sigma_k, \tau_{k+1})}(s)ds,
\]
\[
A(t)=A_{\tau, \sigma}(t):=L^{-1}(t):=\inf\{s\ge0| L(s)\ge t\}.
\]
\begin{defn} We will say that the function
\[
f^{\tau, \sigma}(t):=f(A_{\tau, \sigma}(t)), \ t\ge0
\]
is obtained from  $f$ by deleting the time
interval $\cup_k[\tau_k,\sigma_k).$
\end{defn}
\begin{remk} Informally, we delete  time intervals  $[\tau_k,\sigma_k)$ from the graph of  $f$, then
the remaining time intervals are moved tightly to the left.
\end{remk}

Denote the modulus of continuity of  $f$ by $\omega_f(\delta)=\omega^T_f(\delta):=\mathop{\sup}\limits_{\begin{subarray}{l}
s, t\in[0, T]\\
|s-t|\le\delta\end{subarray}} \rho(f(s),f(t))$.

\begin{lem}\label{lem:lem1}
Assume that
\[
(T+1)-L(T+1)=\int^{T+1}_0\1_{\cup_k[\tau_k, \sigma_k)}(s)ds\le\delta\le1.
\]
Then
\[
\sup_{t\in[0,T]}\rho(f_t,f_t^{\tau, \sigma})\le\omega^{T+1}_f(\delta).
\]
\end{lem}

The proof follows from the fact that  $|A(t)-t|\le\delta, \ t\in[0, T],$ under our assumptions.

\begin{remk} Below we will consider only continuous processes
or c\`{a}dl\`{a}g  processes.
\end{remk}

\begin{thm}\label{lem:lem2}
Let $T>0$. Assume that a sequence of stochastic processes
  $\{X_n, n\ge1\}$ satisfies the following condition
\[
\forall \ \ve>0 \ \exists \ \delta>0 \ \exists \ n_0 \ \forall \ n\ge
n_0
\]
\begin{equation}\label{eq2.2}
\Pb(\omega^{T+1}_{X_n}(\delta)\ge\ve)\le\ve.
\end{equation}
Suppose that for any  $\alpha\in(0,1)$ there exist   sequences of random variables
\[
\tau^{(n)}=\tau^{(n, \alpha)}=\{\tau^{(n, \alpha)}_k, k\ge0\}, \
\sigma^{(n)}=\sigma^{(n, \alpha)}=\{\sigma^{(n, \alpha)}_k, k\ge0\},
\]
that satisfy \eqref{eq1.1}, \eqref{eq1.2} for all $n\ge1$, and there exists a process
$X^{(\alpha)}$ such that
\begin{equation}\label{eq2.3}
\Pb(T+1-L_n^{(\alpha)}(T+1)\ge\alpha)\le\alpha, \ n\ge n_0(T, \alpha);
\end{equation}
\begin{equation}\label{eq2.3a}
X^{(\alpha)}_n\Rightarrow X^{(\alpha)} \ \mbox{as} \ n\to\infty, \text{ in } D([0,T]),
\end{equation}
where
\[
L_n^{(\alpha)}(t)=\int^t_0\1_{\cup_k[\sigma^{(n,\alpha)}_k,\tau^{(n,\alpha)}_{k+1})}(s)\de s,
\ A_n^{(\alpha)}(t)=(L_n^{(\alpha)})^{-1}(t),
\]
and $X^{(\alpha)}_n(t)=X_n(A_n^{(\alpha)}(t))$ is obtained by deleting a time interval from $X_n$.

Then   distributions of $\{X_n, n\ge1\}$ converge weakly in  $D([0,T])$ as $n\to\infty.$
Moreover, distributions of  $\{X^{(\alpha)}\}$ converge weakly in $D([0,T])$ as $\alpha\to0+$,
and the distribution of limits for $\{X^{(\alpha)}\}$ and $\{X_n\}$ coincide.
\end{thm}
\begin{remk} It follows from \eqref{eq2.2} that the limit process has a continuous modification.
\end{remk}
\begin{proof}
It follows from the Skorokhod's theorem  \cite{Skorokhod61-itsp} that for any  $\alpha>0$ there exist copies of random elements
$X^{(\alpha)}$ and $X^{(\alpha)}_n, n\ge1, $ given on the joint probability space $(\wt{\Omega}, \wt{\cF}, \wt{\Pb})$  such that
\begin{equation}\label{eq5.11}
\wt{X}^{(\alpha)}_n\overset{\Pb}{\rightarrow}\wt{X}^{(\alpha)} \text{ in }  \ n\to\infty  \text{ in } D([0,T]).
\end{equation}
We can  construct  a sequence of  random elements
$\{\wt{X}_n, \wt{\tau}^{(n, \alpha)}_k, \wt{\sigma}^{(n, \alpha)}_k, k\ge0\}$
on some extension of $(\wt{\Omega}, \wt{\cF}, \wt{\Pb})$
such that
\[
\{\wt{X}^{(\alpha)}_n, \wt{X}_n, \wt{\tau}^{(n, \alpha)}_k, \wt{\sigma}^{(n, \alpha)}_k, k\ge0\}\overset{d}{=}
\{{X}^{(\alpha)}_n, {X}_n, {\tau}^{(n, \alpha)}_k, {\sigma}^{(n, \alpha)}_k, k\ge0\}
\]
(see reasoning, for example, in  \cite{Kallenberg97}, Ch.5).

Note that $X^{(\alpha)}_n$ can be considered as a measurable function
of  $\{{X}_n, {\tau}^{(n, \alpha)}_k, {\sigma}^{(n, \alpha)}_k, k\ge0\}.$
Thus $\wt{X}^{(\alpha)}_n$ equals a.s. the value 
of the mentioned function
of  $\{\wt{X}_n, \wt{\tau}^{(n, \alpha)}_k, \wt{\sigma}^{(n, \alpha)}_k, k\ge0\}.$
I.e., we obtain $\wt{X}^{(\alpha)}_n$ by deleting   time intervals
$\cup_k[\wt{\tau}^{(n,\alpha)}_k, \wt{\sigma}^{(n,\alpha)}_k)$ from $\wt{X}_n.$

Let $d_0^T$ be the metric in $D([0,T])$ (see for example \cite{Billingsley77}, \S 14).
Observe that  $d_0^T$ is dominated by the uniform metric  on $[0,T]$.

Then Lemma \ref{lem:lem1} yields
\begin{equation}\label{eq6.1}
\begin{split}
&\Pb(d^T_0(\wt{X}_n, \wt{X}^{(\alpha)})\ge2\ve)\le
\\
&\le \Pb(d^T_0(\wt{X}_n, \wt{X}^{(\alpha)}_n)\ge\ve)+\Pb(d^T_0(\wt{X}^{(\alpha)}_n, \wt{X}^{(\alpha)})\ge\ve)\le
\\
&\le \Pb(\sup_{t\in[0,T]}\rho(\wt{X}_n,\wt{X}^{(\alpha)}_n)\ge\ve)+\Pb(d^T_0(\wt{X}^{(\alpha)}_n, \wt{X}^{(\alpha)})\ge\ve)\le
\\
&\le \Pb(T+1-\wt{L}_n^{(\alpha)}(T+1)\ge\alpha)+\Pb(\omega^{T+1}_{\wt{X}_n}(\alpha)\ge\ve)
+\Pb(d^T_0(\wt{X}^{(\alpha)}_n, \wt{X}^{(\alpha)})\ge\ve).
\end{split}
\end{equation}

For a given $\ve>0$ let us select  $\alpha\in(0,\ve)$ and $n_0$ such that (see \eqref{eq2.2} and \eqref{eq5.11})
the second and the third expressions on the right hand side of \eqref{eq6.1}  do not exceed  $\ve$ for $n\ge n_0.$  
Then \eqref{eq2.3} implies that   \eqref{eq6.1} is less than or equal to  $3\ve.$

So, for any $n, m\ge n_0$
\[
\Pb(d^T_0(\wt{X}_n, \wt{X}_m)\ge4\ve)\le \Pb(d^T_0(\wt{X}_n, \wt{X}^{(\alpha)})\ge2\ve)
+\Pb(d^T_0(\wt{X}_m, \wt{X}^{(\alpha)})\ge2\ve)\le6\ve.
\]
Therefore the sequence of distributions of  $\{X_n, n\ge1\}$ is a Cauchy sequence in the Levy-Prokhorov metric (see \cite{EK86}) and hence it is convergent.

Convergence of $\{X^{(\alpha)}, \alpha>0\}$ to the same limit as  $\alpha\to0+$
can be verified  similarly  using bound  \eqref{eq6.1}. Theorem \ref{lem:lem2} is proved.
\end{proof}
\begin{remk} The theorem can be used in the following case. Let
$\{X_n\}$ be a sequence of homogeneous strong Markov processes. Denote by 
 $\tau$   a hitting time of a fixed point  $x_0.$ Assume that the sequence $\{X_n(\cdot\wedge\tau)\}$ converges in distribution to $\{X(\cdot\wedge\tau)\}$ if initial values of the processes converge.
 Let 
$\tau^{(n, \alpha)}_k, \sigma^{(n, \alpha)}_k$ be  sequential entrances into the ball $B(x_0,\alpha/2)$ and exits from the ball $B(x_0,\alpha)$ respectively. Condition 
\eqref{eq2.3} means that if  $\alpha$ is small, then   ``deleted time'' is small uniformly in  $n$. For example, this is true if the average time spent in
$B(x_0,\alpha)$ is small as $\alpha\to 0+$ uniformly in  $n$:
\begin{equation}\label{eqMalVyrezan}
\forall T>0:\ \ \lim_{\alpha\to0+}\sup_n \E\int^T_0\1_{\rho(X_n(t), x_0)\le\alpha}dt=0.
\end{equation}
In some cases the verification of   \eqref{eqMalVyrezan} may be relatively simple  (see, for example, construction of a reflecting Brownian motion in a cone \cite{VW, KwonWilliams}).

Condition   \eqref{eq2.2} about modulus of continuity is similar to a condition that ensures weak relative compactness
of processes in a space of continuous functions. It is not difficult to show that if  $X$ is continuous and
 $X_n(\cdot\wedge\tau)\Rightarrow X(\cdot\wedge\tau)$ whenever initial conditions converge, then
 \eqref{eq2.2} is satisfied if, for example,
 condition \eqref{eqMalVyrezan} holds.
\end{remk}

\section{Convergence of processes that are obtained by deleting   time intervals}

In this section   sufficient conditions that ensure convergence  \eqref{eq2.3a} are given.
Theorem  \ref{prop1} given at the end of this section is an analogue of Theorem  \ref{lem:lem2} that takes into account  these sufficient conditions.

We assume that hitting times of some points are renewal moments  for processes $\{X_n\}$. In this case a distribution of a process with deleted time can be obtained by construction given below.

Consider a space
$$
M=\left\{\ov{t}=\{t_i, i\ge1\}\in (0;\infty)^\mbN:\ \sum_it_i=+\infty\right\}
$$
with the metric of coordinate-wise convergence  $\rho_M(\ov{t},\ov{s})=\sum_n2^{-n}|{t_n}-{s_n}|\wedge1.$

Let $E$ be a locally compact space. Let  $F$ be a map of $C_E([0, \infty))^\mbN  \times M $ to  $D_E([0;\infty))$ defined as
$$
(\ov{f},\ov{t})=(f_1,f_2,\dots;t_1,t_2,\dots)\to F(\ov{f},\ov{t})=\sum_k f_k(\cdot+\sum_{i=1}^{k-1} t_i)
\1_{[\sum_{i=1}^{k-1} t_i,\sum_{i=1}^{k}t_i)}.
$$
Here metric in  $C_E([0, \infty))$  corresponds to the uniform convergence on compact sets, metric on  $C_E([0, \infty))^\mbN$ is introduced similarly to $\rho_M$.

It is not difficult to see that the function $F$ is continuous on   $C_E([0, \infty))^\mbN  \times M $. Therefore, if a sequence of  $C_E([0, \infty))^\mbN  \times M$-valued  random elements  $\{(\ov{\xi}^{(n)},\ov{\tau}^{(n)}), n\geq 0\}$ is such that 
$$
(\ov{\xi}^{(n)},\ov{\tau}^{(n)}) \Rightarrow (\ov{\xi}^{(0)},\ov{\tau}^{(0)}) \text{ as } n\to \infty,
$$
then
$$
F(\ov{\xi}^{(n)},\ov{\tau}^{(n)}) \Rightarrow F(\ov{\xi}^{(0)},\ov{\tau}^{(0)})  \text{ as }  n\to \infty 
$$
in $D_E([0;\infty))$.

\begin{remk}
Everything above  (and also results below) is true if   $\tau_k\in(0,\infty]$ are extended random variables. The corresponding definitions should be  naturally corrected in this case. For example,  if  $\tau_k=+\infty,$ then we should exclude summands with numbers greater than $k$ in the definition of $F(\ov f,\ov t)$.
\end{remk}

 Let $X$ be $E$-valued homogeneous strong Markov process. By $Q_x$ denote the distribution of  $X$ given  $X(0)=x$.
Let $0=\sigma_0\le\tau_1< \sigma_1<\tau_2<\dots$ be a sequence of stopping times such that 
$\sum_k(\tau_{k+1}-\sigma_k)=+\infty$ a.s.
 Set
 $$
 \eta_k=\left\{%
\begin{array}{ll}%
\tau_{k+1}-\sigma_k, & \mbox{ if } \tau_1>0, \\%
\tau_{k+2}-\sigma_{k+1}, & \mbox{ if }\tau_1=0,\\%
\end{array}%
\right.
$$
$$
\xi_k(t)=\left\{%
\begin{array}{ll}%
X((\sigma_k+t)\wedge \tau_{k+1}), & \mbox{ if } \tau_1>0, \\%
X((\sigma_{k+1}+t)\wedge \tau_{k+2}), & \mbox{ if }\tau_1=0. \\%
\end{array}%
\right.
$$

By $X^{(\tau,\sigma)}$ denote a process obtained from  $X$ by deleting time intervals $\cup_k[\tau_k,\sigma_k).$
It can be seen that $X^{(\tau,\sigma)}=F(\ov{\xi},\ov{\eta})$.

Strong Markov property for  $X$ yields that a sequence of random elements 
 $\{(\xi_k,\eta_k), \ k\geq 1\}$ is $C_E([0, \infty))\times(0,\infty)$-valued Markov chain  (may be non-homogeneous). So, we need a result on convergence in distribution for a sequence of Markov chains  (In fact, we will consider not a general  sequence of stopping times  $\{\tau_k,\sigma_k\},$ we will select it in a special way).

Let   $\{X_n\}_{n\ge 1}$ be a sequence of continuous homogeneous strong Markov processes taking values in   $E$. Denote by $Q_x^{(n)}$  the distribution of $X_n$ given $X_n(0)=x$.
Assume that there exists a point  $x^*$ such that for any $n$ the process $X_n$
infinitely often enters any ball $B(x^*, \ve)$ and infinitely often exits from any ball  $B(x^*, \ve)$  with $Q^{(n)}_x$ probability  1 for any
$x\in E.$
Let $\alpha>0$, $\alpha_1\in(0,\alpha)$. By
$\{\tau_k, \sigma_k\}_{k\ge1}$  denote sequential entrances into the ball $B(x^*, \alpha_1)$ and exits from the ball  $B(x^*, \alpha)$, respectively. Set $\sigma^n_0 := 0$,
\begin{equation}\label{eq:tau-sigma}
\begin{split}
\tau^n_{k} = \tau_k^{(n,\alpha_1, \alpha)} := \inf\{t\ge\sigma_{k-1}^{(n,\alpha_1, \alpha)} &\colon \rho(X_n(t),x^*) \le \alpha_1\},\ k\ge 1,
\\
\sigma^n_k = \sigma_k^{(n,\alpha_1, \alpha)} := \inf\{t\ge\tau_k^{(n,\alpha_1, \alpha)} &\colon \rho(X_n(t),x^*) \ge \alpha\},\ k\ge 1.
\end{split}
\end{equation}
Assume that for any  $n\geq 1$ and any initial distribution the following condition is satisfied
\begin{equation}\label{eq:tau-sigma=infty}
\sum_k(\tau^n_{k+1}-\sigma^n_{k})=\infty \ \  \mbox{a.s.}
\end{equation}
Similarly to the previous reasoning let us construct sequences
 $$
(\ov\xi^n,\ov\eta^n)=\{\xi^n_1,\xi^n_2,\dots;\eta^n_1,\eta^n_2,\dots\}= \{(\xi^n_k,\eta^n_k), k\geq 1\}.
$$
 Now the sequence $\{(\xi^n_k,\eta^n_k), k\geq 1\}$ is a homogeneous $C_E([0,\infty))\times(0,\infty)$-valued Markov chain.  Note that its transition probability  $P_n((f,t), A)$ depends only on $f(t)$. Moreover,
\begin{equation}\label{eq:TransKer}
P_n((f,t), A)
=\int_E Q^{(n)}_u\Big((X_n(\cdot\wedge \tau^n_1),\tau^n_1)\in A \Big)
Q^{(n)}_{f(t)}\left(X_n(\sigma^n_1)\in du \right).
\end{equation}
Modifying in a minor way   proofs in  \cite{Karr}, relation  \eqref{eq:TransKer} yields the following result.
\begin{lem}\label{lem:lem_31}
Let  $ \tilde P_0(x,\tilde A), x\in E, \tilde A\in \cB(E)$ and $\bar P_0(x, \bar A), x\in E, \bar A\in \cB(C_E([0,\infty))\times [0,\infty)) $ be stochastic kernels that are continuous in  $x$ (topology of   weak convergence is considered on the space of measures). Assume that the initial distributions $X_n(0)$ of Markov chains converge weakly to a measure $\mu_0$. Suppose also that for any  $x\in E$ and any sequence $\{x_n\}$ that converges   to $x$ we have
\newline
{\bf A1.} $Q^{(n)}_{x_n}\left(X_n(\sigma^n_1)\in \cdot \right)\Rightarrow \tilde P_0(x,\cdot) \text{  as  } n\to\infty,$
\newline
{\bf A2.} $Q^{(n)}_{x_n}\Big((X_n(\cdot\wedge \tau^n_1),\tau^n_1)\in \cdot \Big) \Rightarrow \bar P_0(x,\cdot) \text{  as  } n\to\infty.$

Then the sequence  of  homogeneous Markov chains $\{(\xi^n_k,\eta^n_k), k\geq 1\}$ converges weakly in  $(C_E([0,\infty))\times (0,\infty))^\mbN$ as    $n\to \infty $ to a homogeneous Markov chain
 $\{(\xi^0_k,\eta^0_k), k\geq 1\}$ with the transition kernel (compare with  \eqref{eq:TransKer})
$$
P_0((f,t), A)
=\int_E \bar P_0(u, A)
 \tilde P_0({f(t)}, du).
$$
In particular, if  $X_0$ is a Markov process such that
\begin{equation}\label{eq:X_0}
\mu_0=^d X_0(0),\
 Q^{(0)}_{x}\left((X_0(\cdot\wedge \tau^0_1),\tau^0_1)\in \cdot \right) = \bar P_0(x,\cdot),
\end{equation}
\begin{equation}\label{eq:X_01}
\ Q^{(0)}_{x}\left(X_0(\sigma^0_1)\in \cdot \right)=\tilde P_0(x,\cdot)
\end{equation}
and  $\{(\xi^0_k,\eta^0_k), k\geq 1\}$ are constructed from $X_0$ in the same way as  $\{(\xi^n_k,\eta^n_k), k\geq 1\}$ constructed from $X_n,$ then
$$
(\ov{\xi}^{(n)},\ov{\tau}^{(n)}) \Rightarrow (\ov{\xi}^{(0)},\ov{\tau}^{(0)}) \text{  as  } n\to \infty.
$$
Moreover, we have a convergence of processes  obtained by deleting of  time intervals:
$$
X_n^{(\tau^n,\sigma^n)}\Rightarrow X_0^{(\tau^0,\sigma^0)},\ n\to \infty
$$
in  $D_E([0;\infty))$.
\end{lem}

Combining Theorem \ref{lem:lem2} and obtained sufficient conditions, we get the following results on convergence of stochastic processes.

\begin{thm}\label{prop1}
Let $\{X_n,\ n\ge0\}$ be a sequence of continuous homogeneous strong Markov processes in a locally compact metric space $E$. 
Assume that for any  $T>0$ condition   \eqref{eq2.2} is satisfied, and   for any  $\alpha>0$ there exists  $\alpha_1\in(0,\alpha)$ such that sequences  $\{(\tau_k^{(n,\alpha_1,\alpha)},\sigma_k^{(n,\alpha_1,\alpha)}),k\ge1\},$
(defined in  \eqref{eq:tau-sigma}) satisfy the following conditions: \\
1) \eqref{eq:tau-sigma=infty};\\
2) \eqref{eq2.3} or \eqref{eqMalVyrezan};
\\
3)  conditions {\bf A1}, {\bf A2} of Lemma  \ref{lem:lem_31}.\\
Then distributions of  $\{X_n\}$ converge weakly as $n\to\infty$ in $C([0, \infty))$.

If additionally  \eqref{eq:X_0}, \eqref{eq:X_01} hold  and also 
\begin{equation}\label{eq:eq535}
\int_0^\infty\1_{\{X_0(t)=x^*\}}dt=0 \mbox{a.e.},
\end{equation}
then $X_n\Rightarrow X_0$ as $n\to\infty$ in $C([0, \infty))$.
\end{thm}
\begin{remk} It follows from  \eqref{eq:eq535} that
 $$
 X^{(\tau^{(0,\alpha_1,\alpha)}, \sigma^{(0,\alpha_1,\alpha)})}_0\Rightarrow X_0 \text{  as  } \alpha\to 0+,
 $$
 in $D([0, \infty)).$ Since all processes  $\{X_n\}$ are continuous, then the weak convergence  in  $D([0,T])$ implies the convergence in   $C([0,T])$  as well.
\end{remk}

\begin{remk} Instead of the strong Markov property in Theorem  \ref{prop1} we may assume only that the moments of entrance into any ball and moments of the exit from any ball were ``renewal'' for the processes. For example, this holds if  $E=\mbR$ and $\{X_n\}$ are continuous semi-Markov processes, see the corresponding definition in \cite{Harlamov}.
\end{remk}

\begin{remk} \label{remk_MarkovChain}
The similar statement is also true for processes constructed from Markov chains with a natural correction of formulation. For example, let  $\{X_0(t), t\ge 0\}$ 
be
a continuous strong Markov process with values in  $\mbR^d$,  $\{X^{(n)}(k), k\ge 0\}$ be   $\mbR^d$-valued Markov chains. Put 
$X_n(\frac{k}{n}):=\frac{1}{\sqrt{n}}X^{(n)}(k)$ and define the process  $X_n(t)$ for all  $t\ge 0$ by linearity or in a step-wise manner. Then analogous versions of Lemma   \ref{lem:lem_31} and Theorem  \ref{prop1} are true with the following modifications.

Sequences   $\tau^{(n)}_k, \sigma^{(n)}_k$ should be defined as follows
\begin{equation}\label{eq:tau-sigma1}
\begin{split}
\tau^{(n)}_{k}=\tau^{(n,\alpha_1,\alpha)}_{k}  := \inf\{t\in\tfrac{1}{{n}}\mbZ_+, t\ge\sigma^{(n,\alpha_1,\alpha)}_{k-1} &\colon |X_n(t)-x^*| \le \alpha_1\},\ k\ge 1,\\
\sigma_k^{(n)}=\sigma_k^{(n,\alpha_1,\alpha)} := \inf\{t\in \tfrac{1}{{n}}\mbZ_+, t\ge\tau^{(n,\alpha_1,\alpha)}_k &\colon |X_n(t)-x^*| \ge \alpha\},\ k\ge 1.
\end{split}
\end{equation}
Moreover, if the phase space for $\{X^{(n)}(k), \ k\geq 1\}$ is  $\mbZ^d$ instead of  $\mbR^d,$ then in Lemma \ref{lem:lem_31} we have to consider only sequences $\{x_n\}, \lim_{n\to\infty} x_n=x_0$ such that $x_n\in\frac{1}{\sqrt{n}}\mbZ^d$ for any $n\geq 1.$
\end{remk}
\begin{remk}
It most likely that  \eqref{eq2.2}  is superfluous in Theorem \ref{prop1}, where continuous strong Markov processes are considered. However,
we cannot omit this condition in the case when processes are generated by Markov chains. Theoretically,   $\{X_n\}$ could jump far away from a ball  $B(x^*,\alpha_1)$ by a single  ``large'' jump.

\end{remk}
\section{Limit behavior of perturbed random walk}

Let $\{X(k),\ k\in\Z_+\}$ be a Markov chain on $\Z$
with transition probabilities  $p_{i,j}$ such that for some fixed $m$:
\[
p_{i,i+1} = p_{i,i-1} = 1/2 \quad\text{for}\quad |i|>m
\]
and
\begin{equation}\label{eq:integrable}
\sum_{j\in\Z} |j|\, p_{i,j} < \infty \quad\text{for}\quad |i|\le
m.
\end{equation}

We will call    $X$   a random walk with a membrane $[-m,m]$. Condition  \eqref{eq:integrable} means
that jumps of  $X$ from  the membrane  are integrable.

Extend   the chain $\big\{X(k),\ k\in\Z_+\big\}$ to all $t\ge0$  by linearity.
Define
\[
X_n(t) := \tfrac{1}{\sqrt{n}} X(nt),\ t\ge0,\
n\in\mbN.
\]
In this section we will  prove that distributions of  $\{X_n\}$ converge weakly. 

Let us introduce auxiliary variables. By
\begin{equation}\label{eq:jumps}
\tau := \inf\{k\ge0 \colon |X(k)|>m\}
\end{equation}
denote the exit time from the membrane.
Let $\xi^{(\pm)}$ be a random variable whose  distribution  equals the distribution of $\big(X(\tau)-m\sign X(\tau)\big)$ given  $X(0)=\pm m$.


\begin{thm}\label{thm:main-theorem}
Assume that all states of the chain  $\big\{X(k),\ k\in\Z_+\big\}$ are connected.
Then a sequence of stochastic processes
$\{X_n\}$ converges in distribution in  $C([0,T])$ to a skew Brownian motion
$W_\gamma$, $W_\gamma(0)=0,$ with a parameter
\begin{equation}\label{eq:gamma}
\gamma = \frac{\E\xi^{(+)} \Pb(\xi^{(-)}>0) + \E\xi^{(-)} \Pb(\xi^{(+)}<0)}
   {\E|\xi^{(+)}| \Pb(\xi^{(-)}>0) + \E|\xi^{(-)}| \Pb(\xi^{(+)}<0)},
\end{equation}
i.e., the limit is a continuous Markov process with the transition density 
\begin{equation*}
p_t(x,y) = \varphi_t (x-y) + \gamma\sign(y)\varphi_t (|x|+|y|),\
x,y \in \mathbb{R},
\end{equation*}
where $\varphi_t(x) = \frac{1}{\sqrt{2\pi t}} e^{-x^2/2t}$ is the density of the normal distribution  $N(0,t)$.
\end{thm}
\begin{remk}
The initial distributions   are $X_n(0)=X(0)/\sqrt{n}.$ It is not difficult to prove the result in a scheme of series, where   $X_n(0)=x_n\in \frac{1}{\sqrt{n}}\mbZ$ and $\lim_{n\to\infty}x_n=x.$  In this case, the limiting skew Brownian motion starts from  $x.$
\end{remk}
\begin{remk}
A limit for the sequence $\{X_n\}$ exists even if we omit the assumption that all states of the chain  $\big\{X(k),\ k\in\Z_+\big\}$ are connected.
Then the limit may be not a skew Brownian motion. For example, it could be a Brownian motion that sticks at zero or a mixture of two reflected Brownian motions that reflects at 0 up and down. The list of all possibilities was discussed in  \cite{PP11}. Consideration of those cases is much simpler than the proof of Theorem \ref{thm:main-theorem} below.
\end{remk}

{\it Proof of Theorem \ref{thm:main-theorem}.} Let us apply Theorem~\ref{prop1} or more precisely, its modification for Markov chains  (see Remark
 \ref{remk_MarkovChain}).

Let $\alpha>0$ and $\alpha_1\in(0,\alpha)$ be arbitrary fixed numbers.
Set $\sigma_0^{(n)}:=0$,
\begin{equation}\label{eq:tau-sigma2}\begin{split}
\tau_{k}^{(n)} = \tau_{k}^{(n,\alpha_1,\alpha)}  :=
&\inf\{t\ge\sigma_{k-1}^{(n)},\ t\in\tfrac{1}{{n}}\mbZ_+ \colon |X_n(t)| \le \alpha_1\}=\\
= &\inf\{t\ge\sigma_{k-1}^{(n)},\ t\in\tfrac{1}{{n}}\mbZ_+ \colon |X(nt)| \le \alpha_1\sqrt{n}\},\\
\sigma_k^{(n)} = \sigma_k^{(n,\alpha_1,\alpha)}  := &\inf\{t\ge\tau_k^{(n)},\ t\in\tfrac{1}{{n}}\mbZ_+  \colon
|X_n(t)| \ge \alpha\}=\\
= &\inf\{t\ge\tau_k^{(n)},\ t\in\tfrac{1}{{n}}\mbZ_+  \colon
|X(nt)| \ge \alpha\sqrt{n}\}.
\end{split}\end{equation}

For any $n\ge1$ consider a process  $X_n^{(\alpha)}$  constructed from $X_n$ by deleting time intervals
$\cup_k[\tau_k^{(n,\alpha_1, \alpha)}, \sigma_k^{(n,\alpha_1, \alpha)})$.

Let $X_0$ from Theorem~\ref{prop1} be a skew Brownian motion $W_\gamma$. Note that a skew Brownian motion  $W_\gamma$ behaves as a Wiener process before hitting 0  (see, e.g.,   \cite{HarrShepp81}). So, conditions {\bf A2} and \eqref{eq:X_0} reformulated for Markov chains are satisfied.
Condition \eqref{eq:eq535} also holds  because a skew Brownian motion has a transition density. 

\subsection{Distribution at an instant of exit time}\label{sec:xtilde-ytilde}
Let us verify conditions   {\bf A1} and \eqref{eq:X_01} of Lemma   \ref{lem:lem_31} for the sequence   $\{X_n, n\ge0\}$, where  $X_0$ equals  $W_\gamma$.
At first, let us find the conditional distribution of  $W_\gamma(\sigma_k^{(0)})$ given $W_\gamma(\tau_k^{(0)})=\alpha_1$.
This distribution is concentrated at $\pm\alpha.$ The corresponding weights are equal to probabilities  for $W_\gamma$ to exit
from 
$[-\alpha,\alpha]$ through the right or left ends, respectively, given    $W_\gamma(0)=\alpha_1$.

It is known that (see, e.g.,   \cite{Lejay06}) the scale function for the skew Brownian motion 
$W_\gamma$ equals
\[
\psi(x) = \left\{%
\begin{array}{ll}%
x/p, & x\ge0, \\%
x/q, & x<0, \\%
\end{array}%
\right.
\]
where $p=\frac{1+\gamma}{2}, q=\frac{1-\gamma}{2}.$
Thus,
\[
\Pb(W_\gamma(\sigma^{(0)}_k) = +\alpha ~|~ W_\gamma(\tau^{(0)}_k) = +\alpha_1) =
\frac{\psi(\alpha_1) - \psi(-\alpha)}{\psi(\alpha) - \psi(-\alpha)},
\]
\[
\Pb(W_\gamma(\sigma^{(0)}_k) = -\alpha ~|~ W_\gamma(\tau^{(0)}_k) = +\alpha_1) =
\frac{\psi(\alpha) - \psi(\alpha_1)}{\psi(\alpha) - \psi(-\alpha)},
\]
If   $W_\gamma(\tau_k^{(0)})=-\alpha_1$, then the corresponding distribution is
\[
\Pb(W_\gamma(\sigma^{(0)}_k) = +\alpha ~|~ W_\gamma(\tau^{(0)}_k) = -\alpha_1) =
\frac{\psi(-\alpha_1) - \psi(-\alpha)}{\psi(\alpha) - \psi(-\alpha)},
\]
\[
\Pb(W_\gamma(\sigma^{(0)}_k) = -\alpha ~|~ W_\gamma(\tau^{(0)}_k) = -\alpha_1) =
\frac{\psi(\alpha) - \psi(-\alpha_1)}{\psi(\alpha) - \psi(-\alpha)}.
\]


Let us find the distribution of $X_n(\sigma_k^{(n)})$.

Set  $C_n = -[-\alpha\sqrt{n}],$ where $[\cdot ]$ denotes the integer part.

Let   $X(0)=i\in\{-C_n,\ldots,+C_n\}$. By $\rho^{(n)}_i$ denote the probability that  $\{X(k)\}$ hits $+C_n$ before it visits the set $(-\infty,C_n]\cup(C_n,+\infty)$.
These probabilities satisfy the system of linear equations:
\[
\left\{%
\begin{array}{l}
    \rho^{(n)}_{C_n} = 1, \\
    \rho^{(n)}_i = \rho^{(n)}_{-C_n} = 0,\ |i|>C_n, \\
    \rho^{(n)}_i = \tfrac{1}{2} \rho^{(n)}_{i-1} + \tfrac{1}{2} \rho^{(n)}_{i+1},\ m<|i|<C_n, \\
    \rho^{(n)}_{\pm m} = \Sum_{m< |j|\le C_n} p_{\pm m, j} \rho^{(n)}_j, \\
\end{array}%
\right.
\]
where $p_{\pm m,j} = \Pb\big(\xi^{(\pm)} = j-m\sign(j)\big)$.

Consider $\rho^{(n)}_i$ for $m\le i\le C_n$.
Note that all points
\[
(m, \rho^{(n)}_{m}),\ (m{+}1, \rho^{(n)}_{m+1}),\ \ldots,\ (C_n, \rho^{(n)}_{C_n})
\]
belong to the same line. Therefore,
\[
\rho^{(n)}_{m+k} = \rho^{(n)}_{m} + \tfrac{k}{C_n'} (1-\rho^{(n)}_{m}) = \rho^{(n)}_{m}
(1-\tfrac{k}{C_n'}) + \tfrac{k}{C_n'},\ k=\overline{0,C_n'},
\]
where $C_n'=C_n-m$.

Similarly,
\[
\rho^{(n)}_{-m-k} = \rho^{(n)}_{-m} (1-\tfrac{k}{C_n'}),\ k=\overline{0,C_n'}.
\]

Substituting these expressions into the equations for  $\rho^{(n)}_{\pm m}$
and then rewriting the obtained systems in terms of $\xi^{(+)}$ and
$\xi^{(-)}$, we get
\begin{multline*}
\rho^{(n)}_{m} \big(C_n \Pb(\wt\xi_n^{(+)} < 0) + \E(\wt\xi_n^{(+)}\vee0)\big) = \\
= \rho^{(n)}_{-m} \big(C_n \Pb(\wt\xi_n^{(+)} < 0) + \E(\wt\xi_n^{(+)}\wedge0)\big) + \E(\wt\xi_n^{(+)}\vee0),
\end{multline*}
\begin{multline*}
\rho^{(n)}_{-m} \big(C_n \Pb(\wt\xi_n^{(-)} > 0) + \E(\wt\xi_n^{(-)}\wedge0)\big) = \\
= \rho^{(n)}_{m} \big(C_n \Pb(\wt\xi_n^{(-)} > 0) - \E(\wt\xi_n^{(-)}\vee0)\big) + \E(\wt\xi_n^{(-)}\vee0),
\end{multline*}
where $\wt\xi_n^{(\pm)} 
:= \xi^{(\pm)}\1_{|\xi^{(\pm)}|\le C_n'}$.
Thus,
\[
\rho^{(n)}_{m} = \frac{C_n \Pb(\wt\xi_n^{(+)}<0) \E(\wt\xi_n^{(-)}\vee0) + C_n
\Pb(\wt\xi_n^{(-)}>0) \E(\wt\xi_n^{(+)}\vee0) + A_n}{C_n \Pb(\wt\xi_n^{(+)}<0) \E|\wt\xi_n^{(-)}| + C_n \Pb(\wt\xi_n^{(-)}>0)
\E|\wt\xi_n^{(+)}| + A_n},
\]
where $A_n = \E(\wt\xi_n^{(+)}\wedge0) \E(\wt\xi_n^{(-)}\vee0) - \E(\wt\xi_n^{(+)}\vee0) \E(\wt\xi_n^{(-)}\wedge0)$.

Analogously,
\[
\rho^{(n)}_{-m} = \frac{C_n \Pb(\wt\xi_n^{(+)}<0) \E(\wt\xi_n^{(-)}\vee0) + C_n
\Pb(\wt\xi_n^{(-)}>0) \E(\wt\xi_n^{(+)}\vee0)}{C_n \Pb(\wt\xi_n^{(+)}<0) \E|\wt\xi_n^{(-)}| + C_n
\Pb(\wt\xi_n^{(-)}>0) \E|\wt\xi_n^{(+)}| + A_n}.
\]

It can be easily checked that for any $\alpha>0$
\[
\Pb\big(\wt \xi_n^{(\pm)} \neq \xi^{(\pm)}\big) \le
\Pb\big(|\xi^{(\pm)}|>C_n]\big) \to 0  \text{ as }\ n\to\infty
\]
and
\begin{multline}\label{eq:lim-pn}
\Lim_{n\to\infty} \rho^{(n)}_{m} = \Lim_{n\to\infty} \rho^{(n)}_{-m} = \\
= p =
\frac{\Pb(\xi^{(+)}<0) \E(\xi^{(-)}\vee0) + \Pb(\xi^{(-)}>0) \E(\xi^{(+)}\vee0)}{\Pb(\xi^{(+)}<0)
\E|\xi^{(-)}| + \Pb(\xi^{(-)}>0) \E|\xi^{(+)}|}.
\end{multline}

The  above reasoning implies that for any  $\alpha>0$ we have the uniform convergence
\[
\sup_{m\le |i|\le C_n} \Big|\rho_i^{(n)}-\frac{\psi(\frac{i}{n}) - \psi(-\alpha)}{\psi(\alpha) - \psi(-\alpha)}\Big| \to 0 \text{ as } \ n\to\infty.
\]

The similar formulas are also true for the probabilities of hitting  $-\alpha.$
In particular, these formulas yield that the probability to exit from   $[-\alpha,\alpha]$
by a large jump from the membrane goes to 0 as  $n\to\infty.$

Thus, conditions  {\bf A1} and \eqref{eq:X_0} are satisfied.

\subsection{Modulus of continuity of  $X_n$}\label{sec:mod-xn}

Let us verify condition \eqref{eq2.2} for processes  $\{X_n\}$. It is sufficient to show that
\[
\forall\ T>0\ \forall \ \ve>0 \ \exists \ \delta>0 \ \exists \ n_0 \ \forall \ n\ge n_0
\]
\[
 \Pb(\omega^{T}_{X_n}(\delta)\ge\ve)\le\ve,
\]
where $\omega_f(\delta) = \omega^T_f(\delta)
$ is the modulus of continuity of a function $f$ on $[0,T]$.

Compare the distribution of the modulus of continuity of a process   $X_n$
with the distribution of the modulus of continuity of a symmetric random walk with unit jumps. To do this, we construct copies of these processes on the same probability space.

Let $\{S(k)\}$ be a symmetric  random walk with unit jumps that is independent of a chain $X$.

Let  $\xi_k$ be jumps of  $X$ from the membrane:
\[
\xi_k = X(\tau_k) - m\sign(X(\tau_k)),\ k\ge1,
\]
where $\tau_k$ is the $k$-th   exit of   $X$ from  $[-m,m]$.

Introduce the following auxiliary process. Set
\[
\wt X(k) := S(k),\ k=\overline{0,t_1},
\]
where $t_1$ is the first instant of hitting 0 for  $\wt X$.\\
Define the next increment for  $\wt X$ by $\xi_1$:
\[
\wt X(t_1+1) := \wt X(t_1) + \xi_1 = 0 + \xi_1 = \xi_1.
\]
Further increments of  $\wt X$ are the same as for $S$ until the instant $t_2$ of the next hitting zero by  $\wt X$:
\[
\wt X(k) := \wt X(t_1+1) + S(k-1),\ k=\overline{t_1+2,t_2}.
\]
Set
\[
\wt X(t_2 + 1) := \xi_2,
\]
and so on.

We have a representation
\[
\wt X(k) = S(k-r(k)) + \sum_{j=1}^{r(k)} \xi_j,
\]
where $r(k)=r_{\wt X}(k)$ is the number of visits of zero by the sequence  $\{\wt X(l),\
l=\overline{0,k}\}$.

Set
\[
S_n(t) = \tfrac{1}{\sqrt{n}} \big(S([nt]) + (nt-[nt])S([nt]+1)\big),\ t\ge0,
\]
\[
\wt X_n(t) = \tfrac{1}{\sqrt{n}} \big(\wt X([nt]) + (nt-[nt])\wt X([nt]+1)\big),\ t\ge0.
\]

It can be seen  that the  distribution of the following processes coincide
\begin{equation}\label{eq:xtilde1}
\wt X_n(t) \peq X_n^{(\tau,\wt\tau)}(t) - \tfrac{m}{\sqrt{n}}\sign(X_n^{(\tau,\wt\tau)}(t)),\ t\ge0,
\end{equation}
where $X_n^{(\tau,\wt\tau)}$ is the process obtained from  $X_n$ by deleting time intervals  $\bigcup_{k\ge1}[\tau_k,\wt\tau_k)$. Here  $\wt\tau_0=0$,
\[
\tau_{k} = \tau_{k}^{(n)}  :=
\inf\{t\ge\wt\tau_{k-1}^{(n)},\ t\in \frac{1}{\sqrt{n}}\mbN \colon |X_n(t)| \le m/\sqrt{n}\},\
k\ge 1,
\]
\[
\wt\tau_{k} = \wt\tau_{k}^{(n)}  :=
\inf\{t\ge\tau_k^{(n)},\ t\in \frac{1}{\sqrt{n}}\mbN \colon |X_n(t)| \ge (m+1)/\sqrt{n}\},\
k\ge 1.
\]

It follows from the construction of a process with deleted time intervals that
\begin{equation}\label{eq:xtilde2}
\omega_{X_n}^{~}(\delta) \le \omega_{X_n^{(\tau,\wt\tau)}}(\delta)+2m/\sqrt{n},\ \delta>0.
\end{equation}

Let us compare modulus of continuity for the sequences    $\{S(k)\}$ and $\{\wt X(k)\}$,
and then for the processes $\{S_n(t)\}$ and $\{\wt X_n(t)\}$.
Consider the difference $\wt X(l) - \wt X(k)$, where $l$ and $k>l$  are integers from  $[0,nT]$,
 $r(n)=r_{\wt X}(n)$ is the number of visits of zero by the sequence  $\{\wt X(i),\ i=\overline{0,n}\}$.

Observe that for any $p\in\mbN$ and $0\le k<l\le nT$
\[
\sup_{|l-k|\le p} |\wt X(l) - \wt X(k)| \le 2 \sup_{|l-k|\le p} |S(l) - S(k)| +
\sup_{j\le r(nT)} |\xi_j|.
\]
Indeed, if there was no  visit of zero by $\{\wt X(j),j\in[k,l]\}$, then 
$|\wt X(l) - \wt X(k)| $ does not exceed   $\sup_{|j-i|\le |l-k|} |S(j) - S(i)|$ by construction.
Otherwise, let $k_1:=\inf\{j\geq k\ :\ \wt X(j)=0\}$, $l_1:=\sup\{j\leq l\ :\ \wt X(j)=0\}$. Then
$$
|\wt X(l) - \wt X(k)|\le |\wt X(k_1) - \wt X(k)|+|\wt X(l_1) - \wt X(k_1)|+|\wt X(l_1+1) - \wt X(l_1)|+
$$
$$
+|\wt X(l) - \wt X(l_1+1)|= |\wt X(k_1) - \wt X(k)|+ |\wt X(l_1+1) - \wt X(l_1)|+|\wt X(l) - \wt X(l_1+1)|\le
$$
$$
 \le 2 \sup_{|j-i|\le |l-k|} |S(j) - S(i)|+
\sup_{j\le r(nT)} |\xi_j|.
$$

So, for any    $\delta>0$ we have
\begin{equation}\label{eq:xtilde-s-xi}
\omega_{\wt X_n}(\delta) \le 2 \omega_{S_n}(\delta) + \tfrac{1}{\sqrt{n}} \sup_{j\le r(nT)} |\xi_j|.
\end{equation}


Therefore, for any $\alpha>0$
\begin{equation}\label{eq:xtilde-s-xi-prob}
\Pb(\omega_{\wt X_n}(\delta) > \alpha) \le \Pb(\omega_{S_n}(\delta) > \alpha/3) +
\Pb(\tfrac{1}{\sqrt{n}} \sup_{j\le r(nT)} |\xi_j| > \alpha/3).
\end{equation}

It follows from the weak convergence of $S_n$ in  $\Ce$ (Donsker's theorem)
and  \cite[Theorem 8.2]{Billingsley77}
that
\[
\forall \varepsilon>0\ \forall \alpha>0\ \exists \delta>0\
\exists n_1 \ \forall n\ge n_1:
\]
\[
\Pb(\omega_{S_n}(\delta) > \alpha/3) < \varepsilon/2.
\]

To estimate the second summand in \eqref{eq:xtilde-s-xi-prob},
let us verify that for any  $\delta>0$
\[
\Pb(\max_{j\le r(n)} |\xi_j| > \delta\sqrt{n}) \to 0 \text{ as } \
n\to\infty.
\]

{
For $x>0$, bound
\begin{equation}\label{eq:estim-max-xi}
\Pb(\max_{j\le r(n)} |\xi_j| > \delta\sqrt{n}) \le
\Pb(r(n) > x\sqrt{n}) + \Pb(\max_{j \le x\sqrt{n}} |\xi_j| > \delta\sqrt{n}).
\end{equation}

To bound the first item on the right hand side of  \eqref{eq:estim-max-xi}, we will study a random walk $\wt{S}$ with unit jumps which is  constructed in a special way using trajectories of $\wt X$. Then we compare the number of visits of zero by 
  $\wt{S}$ and $X$, see below. Informally, if we have a greater jump from  zero, then we need more time to return to zero again. So the number of visits 
of zero would be less. 

The formal construction is the following.
Define
\[
\wt S(k) := \wt X(k),\ k=\overline{0,\wt t_1},
\]
where $\wt t_1 := \inf\{k\colon \wt S(k)=0\}$ is the  first instant of   hitting zero by  $\wt S$.
Furthermore, define
\[
\wt S(\wt t_1+1) := \sign(\wt X(t_1+1)),
\]
where $t_1$ is the  first instant of   hitting zero by $\wt X$.

Let increments of $\wt S$ be the same as  the corresponding increments of $\wt X$ until $\wt t_2$, where $\wt t_2$ is the second
hitting time of zero by  $\wt S$:
\[
\wt S(\wt t_1+1+k) := \wt X(t_1+1+k) - \wt X(t_1+1),\ k=\overline{1,\wt t_2 {-} \wt t_1 {-} 1},
\]
and
\[
\wt S(\wt t_2+1) := \sign(\wt X(t_2+1)).
\]
Further construction can be done in the similar way.

Since
$|\wt X(t_i+1)| \ge 1=|\wt S(\wt t_i+1)|$ and
$|\wt X(k+1)-\wt X(k)|=1,\ k\notin\{t_i\}$,
we have $\wt t_{i+1}-\wt t_{i} \le t_{i+1}-t_i,\ i\ge1$.
So, the number of visits of zero  
by $\wt X$ is less than or equal to the number of visits of zero by $\wt S$:
\begin{equation}\label{eq:X_S}
r_{\wt{X}} (k) \le r_{\wt{S}}(k),\ k \ge 0.
\end{equation}

Note that  $\{|\wt S(k)|,\ k=\overline{0,n}\},$ has the same distribution as the absolute value of a symmetric random walk. Therefore, the number of visits of zero
for $\wt S$
(or $|\wt S|$) coincides with the one for a symmetric random walk  $S$ with unit jumps.
The asymptotics for this is well known  (see, e.g.,  \cite{Feller67}). Hence,
\begin{equation}\label{eq:r-of-n-asymp}
\forall x>0  \   \lim_{n\to\infty}\Pb(r_{\wt S}(n)\le x\sqrt{n})=
\sqrt{\tfrac{2}{\pi}} \int_0^x e^{-z^2/2} \de z.
\end{equation}
It follows from \eqref{eq:X_S} and \eqref{eq:r-of-n-asymp} that\\

$\forall \varepsilon>0$ $\exists x>0$ $\exists n_0
$ $\forall n\ge n_0$:
\[
\Pb(r_{\wt X}(n) > x\sqrt{n}) \le \Pb(r_{\wt S}^{~}(n) > x\sqrt{n}) < \varepsilon.
\]

We have as $n\to\infty$
\begin{multline*}
\Pb\big(\max_{j\le x\sqrt{n}} |\xi_j| > \delta\sqrt{n}\big) =
\Pb\big(\bigcup_{j\le x\sqrt{n}} \{|\xi_j| > \delta\sqrt{n} \}\big)
\le \sum_{j \le x\sqrt{n}} \Pb\big(|\xi_j| > \delta\sqrt{n}\big) \le\\
\le
[x\sqrt{n}]\left\{ \Pb\left(|\xi^{(+)}| > \delta\sqrt{n}\right)+\Pb\left(|\xi^{(-)}| > \delta\sqrt{n}\right)\right\}\to 0, 
\end{multline*}
because $\E\xi^{(\pm)}<\infty$ by assumption of the Theorem.


}

It follows from  \eqref{eq:estim-max-xi} and \eqref{eq:xtilde-s-xi-prob} that
\[
\forall \varepsilon>0\ \forall \alpha>0\
\exists \delta>0\ \exists n_0 \ \forall
n\ge n_0:
\]
\begin{equation}\label{eq:zvezdochka}
\Pb\big(\omega_{\wt X_n}(\delta) > \alpha\big) < \varepsilon.
\end{equation}

We have (see \eqref{eq:xtilde1})
\[
\Pb\big(\omega_{X_n^{(\tau,\wt\tau)}}(\delta) > \alpha\big) \le
\Pb\big(\omega_{\wt X_n}(\delta) > \alpha - 2m/\sqrt{n}\big) \le
\Pb\big(\omega_{\wt X_n}(\delta) > \alpha/2\big)
\]
for $n>16m^2/\alpha^2$. So
  \eqref{eq:xtilde2} yields
\[
\forall \ \ve>0 \ \exists \ \delta>0 \ \exists \ n_0 \
\forall \ n\ge n_0:
\]
\[
 \Pb(\omega_{X_n}(\delta)\ge\ve) \le \Pb(\omega_{X_n^{(\tau,\wt\tau)}}(\delta)\ge\ve/2) \le \ve.
\]
Thus, condition \eqref{eq2.2} holds.

\subsection{Investigation of a time spent in a neighborhood of zero}

In this Section we prove that
\begin{equation}\label{eq_fin5}
\forall T>0\ \forall\delta>0\ \  \lim_{\alpha\to 0+}\lim_{n\to\infty}\Pb\left(\int_0^T\1_{|X_n(t)|\le\alpha}dt>\delta\right)=0.
\end{equation}
This, in particular, implies that the sequence of processes  $\{X_n, n\ge 1\}$ satisfies  \eqref{eq2.3}.

Similarly to the reasoning of Section  \ref{sec:mod-xn}, construct processes $\wt S_n,  \wt X_n, X_n^{(\tau,\wt\tau)}.$
The time spent by 
$\{\wt S_n(t), t\in[0,T]\}$ in $[-\alpha, \alpha]$
does not exceed the similar time for the process $\wt X_n.$
So, for all  $T>0, \delta>0$
\begin{equation}\label{eq:that}
  \lim_{\alpha\to 0+}\lim_{n\to\infty}\Pb\left(\int_0^T\1_{|\wt X_n(t)|\le\alpha}dt>\delta\right)\le \lim_{\alpha\to 0+}\lim_{n\to\infty}\Pb\left(\int_0^T\1_{|\wt S_n(t)|\le\alpha}dt>\delta\right).
\end{equation}
Since $\{|\wt S_n(t)|, t\in[0,T]\}$ converges in distribution to a reflected Brownian motion, the right hand side of  \eqref{eq:that} equals $0$.
Therefore,
$$
\lim_{\alpha\to 0+}\lim_{n\to\infty}\Pb\left(\int_0^T\1_{|X_n(t)|\in(\frac{m}{\sqrt{n}},\alpha]}dt>\delta\right)\le
\lim_{\alpha\to 0+}\lim_{n\to\infty}\Pb\left(\int_0^T\1_{|\wt X_n(t)|\le\alpha}dt>\delta\right)=0.
$$

Thus, to prove   \eqref{eq_fin5}, and hence to prove Theorem  \ref{thm:main-theorem}, it suffices to verify
\begin{equation}\label{eq_fin51}
\forall T>0\ \forall\delta>0: \  \lim_{n\to\infty}\Pb\left(\int_0^T\1_{|X_n(t)|\le \frac{m}{\sqrt{n}}}dt>\delta\right)=0.
\end{equation}
Assume that random variables $\zeta^{(+)}$ and $\zeta^{(-)}$ have the same distribution as the time spent by  $X$  in the membrane $\{-m,-m+1,\dots,m\}$ if $X$ enters this set through  $m$ or $-m$, respectively.
By $r_X(k)$ denote the number of entrances of $\{X(i), 0\le i\le k\}$ in the membrane. Let
 $\zeta^{(+)}_j$, $\zeta^{(-)}_j,\ j\ge 0,$ be independent copies of  $\zeta^{(\pm)}$ that are also independent of
 $\wt S.$ It can be easily seen that for any $x>0,\ k\in\mbN:$
 $$
 \Pb(r_X(k)>x)\le \Pb(r_{\wt S}(k)>x),
 $$
 $$
 \Pb\left(\sum_{i=1}^{r_X(k)} \1_{|X(i)|\le m}>x\right)\le
 \Pb\left(\sum_{i=1}^{r_{\wt S}(k)}(\zeta_i^{(+)}+\zeta_i^{(-)})>x\right).
 $$
Therefore
 $$
 \lim_{n\to\infty}\Pb\left(\int_0^T\1_{|X_n(t)|\le \frac{m}{\sqrt{n}}}dt>\delta\right)\le
 \lim_{n\to\infty} \Pb\left(\sum_{i=1}^{r_{\wt S}([nT]+1)}\frac{\zeta_i^{(+)}+\zeta_i^{(-)}}{n}>\delta\right).
 $$
 Arguing as in \S \ref{sec:mod-xn}, the last inequality and   \eqref{eq:r-of-n-asymp} imply  \eqref{eq_fin51}. Theorem  \ref{thm:main-theorem} is proved.

\end{document}